\documentclass{article}
\usepackage{times}

\usepackage{latexsym,amssymb,lastpage}
\usepackage{graphicx,amsfonts}
\usepackage{graphicx}
\usepackage{mathptmx,bm,amsmath}

\newtheorem{theorem}{\bf Theorem}[section]

\numberwithin{equation}{section}

\title{Numerical Solutions of \\
       Matrix Differential Models using \\
       Cubic Matrix Splines II%
}
\date{}
\author{E.\ Defez\footnote{Author to whom all correspondence should
        be addressed.}, A.\ Herv\'{a}s, L.\ Soler, M.M.\ Tung\\
   \footnotesize Instituto de Matem\'{a}tica Multidisciplinar \\
   \footnotesize Universidad Polit\'{e}cnica de Valencia, Spain \\
   \footnotesize $\left\{ \right.$edefez, ahervas, mtung$\left. \right\}$@imm.upv.es}
\begin{document}
\maketitle

\begin{abstract}
\noindent
This paper presents the non-linear generalization of a previous
work on matrix differential models~\cite{mm}. It focusses on the
construction of approximate solutions of first-order matrix
differential equations $Y'(x)=f(x,Y(x))$ using matrix-cubic splines.
An estimation of the approximation error, an algorithm for its
implementation and illustrative examples for Sylvester and Riccati
matrix differential equations are given.
\end{abstract}

\vspace{3mm}


\vspace{2mm}

\section{Introduction}\label{seccion1}

A great variety of phenomena in physics and engineering can be modelled
in the form of matrix-differential equations. Although linear
matrix-differential equations, whose numerical solutions using cubic matrix
splines were presented in~\cite{mm}, are valid for a wide range of applications,
non-linear equations are also of great interest. This work generalizes the
approach of~\cite{mm}, providing a novel scheme to numerically solve non-linear
differential matrix equations of the first-order.
Concretely, in this work we will develop a method for the numerical
integration of the first order matrix differential equation
given by



\begin{equation}\label{problem}
\left.  \begin{array}{rcl}
Y'(x) & = & f(x, Y(x))   \\
\\
Y(a)  & = & Y_{a}
\end{array}
 \right\}  \ a \leq x \leq b \ ,
\end{equation}

where $Y_{a}, Y(t) \in {\mathbb C}^{r \times q}$,
$f:[a,b] \times {\mathbb C}^{r \times q} \mapsto {\mathbb C}^{r \times q}$.\\

Different examples of problem (\ref{problem}) can be found in
\cite{21}. Numerical schemes to obtain approximate solutions
for~(\ref{problem}) by means of linear multistep methods with
constant steps have been devised in~\cite{15b}. Although there exist
{\it a priori} error bounds for these methods expressed in function
of the data problem, these error bounds are given in terms of an
exponential which depends on the integration step $h$. Therefore, in
practice, $h$ will take too small values. Furthermore, these methods
require some interpolation techniques in
order to get a continuous solution \cite{15b}.\\

Generalizing the method proposed for the linear case in~\cite{mm},
here we elaborate an extension using cubic-matrix splines in the
numerical approximation for the solutions of (\ref{problem}). In the
scalar case, cubic splines were used in \cite{20} for the resolution
of ordinary differential equations obtaining approximations that,
among other advantages, were of class ${\cal C}^1$ in the interval
$[a,b]$. These splines are easy to compute and produce an
approximation error of only $O(h^4)$. Recently, this method has been
used in the resolution of other scalar problems as discussed in
\cite{nueva6}, and even linear matrix problems (see \cite{mm}). The
present work extends this powerful scheme to
the resolution of matrix problems of the non-linear type (\ref{problem}).\\

This paper is organized as follows. In section \ref{seccion2} we
develop the proposed method, whose algorithm is then given in
Section \ref{seccion3}. Finally, in Sections \ref{seccion7},
\ref{seccion4} and \ref{seccion5} practical examples are
presented.\\

Throughout this work, we will adopt the notation for norms and
matrix cubic splines as in the previous work~\cite{mm} and common in
matrix calculus. Following this nomenclature, we define the
Kronecker product of $A \ = \ \left(a_{ij} \right) \in {\mathbb
C}^{m \times n}$ and $B \in {\mathbb C}^{r \times s}$, denoted by $A
\otimes B$, as the block matrix
$$
A \otimes B \ = \ \left[\begin{array}{ccc}
a_{11} B & \ldots & a_{1n} B \\
\vdots  &         & \vdots \\
a_{m1} B & \ldots & a_{mn} B
 \end{array} \right] \ .
$$
The column-vector operator on a matrix $A \in {\mathbb C}^{m \times
n}$ is given by

$$
vec(A) \ = \ \left[\begin{array}{c}
A_{\bullet 1} \\
\vdots \\
A_{\bullet n}\end{array} \right] \ , \ \mbox{where} \
A_{\bullet k} \ = \ \left[\begin{array}{c}
a_{1k} \\
\vdots \\
a_{mk}\end{array} \right] \ .
$$

If $Y \ = \ \left(y_{ij} \right) \in {\mathbb C}^{p \times q}$ and $X \
= \ \left(x_{ij} \right) \in {\mathbb C}^{m \times n}$, then the
derivative of a matrix with respect to a matrix is defined by
\cite[p.62 and 81]{graham}:
$$
\frac{\partial Y}{\partial X} \ = \ \left[\begin{array}{ccc}
\displaystyle \frac{\partial Y}{\partial x_{11}} & \ldots &\displaystyle  \frac{\partial Y}{\partial x_{1n}} \\
\vdots  &         & \vdots \\
\displaystyle \frac{\partial Y}{\partial x_{m1}} & \ldots &\displaystyle  \frac{\partial Y}{\partial x_{mn}}
\end{array} \right] \ ,\ \mbox{where}\quad
\ \ \frac{\partial Y}{\partial x_{rs}} \ = \ \left[\begin{array}{ccc}
\displaystyle \frac{\partial y_{11}}{\partial x_{rs}} & \ldots & \displaystyle  \frac{\partial y_{1q}}{\partial x_{rs}} \\
\vdots  &         & \vdots \\
\displaystyle \frac{\partial y_{p1}}{\partial x_{rs}} & \ldots &\displaystyle  \frac{\partial y_{pq}}{\partial x_{rs}}
 \end{array} \right] \ . \
$$
If $X \in {\mathbb C}^{m \times n}, Y \in {\mathbb C}^{n \times v}, Z
\in {\mathbb C}^{p \times q}$, then the following rule for the
derivative of a matrix product with respect to another matrix applies
\cite[p.84]{graham}:

\begin{equation}\label{tema0}
\frac{\partial XY}{\partial Z} \ = \
\frac{\partial X}{\partial Z}\left[I_{q} \otimes Y\right] \ + \
\left[I_{p} \otimes X \right] \frac{\partial Y}{\partial Z} \ ,
\end{equation}

where $I_{q}$ and $I_{p}$ denote the identity matrices of dimensions
$q$ and $p$, respectively. If $X \in {\mathbb C}^{m \times n}, Y \in
{\mathbb C}^{u \times v}, Z \in {\mathbb C}^{p \times q}$, the
following chain rule \cite[p.88]{graham} is valid :
\begin{equation}\label{tema1}
\frac{\partial Z}{\partial X} \ = \
\left[\frac{\partial \left[vec(Y)
\right]^{t}}{\partial X} \otimes I_{p}\right]
\left[I_{n} \otimes \frac{\partial Z}{\partial \left[vec(Y)
\right]}\right] \ .
\end{equation}

If $A \ = \ \left(a_{ij} \right) \in {\mathbb C}^{m \times n}$, the
Frobenius norm of $A$ is \cite{4} given by:

\begin{equation}\label{frobenius}
\left\|A \right\|_F=\sqrt{\sum_{i=1}^{m}\sum_{j=1}^{n} \left|a_{ij}
\right|^{2}} \ .
\end{equation}

The following relationship between the 2-norm and Frobenius norm
holds \cite{4}:

\begin{equation}\label{frobenius1}
\left\|A \right\|_2 \leq \left\|A \right\|_F \leq \sqrt{n} \left\|A
\right\|_2 \ .
\end{equation}

\section{Proposed general method}\label{seccion2}
Let us consider the problem

\begin{equation}\label{tema2}
\left.  \begin{array}{rcl}
Y'(x) & = & f(x, Y(x))   \\
Y(a)  & = & Y_{a}
\end{array}
 \right\} \ a \leq x \leq b \ ,
\end{equation}
where $Y_{a}, Y(t) \in {\mathbb C}^{r \times q}$,
$f:[a,b] \times {\mathbb C}^{r \times q} \mapsto {\mathbb C}^{r \times q}$,
$f \in {\cal C}^{1}\left(T \right)$, with
\begin{equation}\label{guapa1}
T \ = \ \left\{(x,Y) \ ; \
a \leq x \leq b \ , \ Y \in {\mathbb C}^{r \times q} \right\} \ ,
\end{equation}
and $f$ fulfills the global Lipschitz's condition
\begin{equation}\label{2}
\left\|f\left(x, Y_{1}\right) \ - \ f\left(x, Y_{2}\right)\right\| \ \leq L
\left\|Y_{1} - Y_{2}   \right\| \ , \  a \leq x \leq b \ ,
Y_{1}, Y_{2} \in {\mathbb C}^{r \times q} \ ,
\end{equation}
which guarantees the existence and uniqueness of the continuously
differentiable solution $Y(x)$ of problem (\ref{tema2}), see \cite[p.99]{flett}.\\

Let us consider $h=(b-a)/n$, $n$ being a positive integer, so that
the partition of the interval $[a,b]$ is given by
\begin{equation}\label{partition}
\Delta_{[a,b]}=\left\{a=x_{0} \ < \ x_{1} \ < \ \ldots \ < \ x_{n}=b
\right\} \ , \ x_{k}=a+ k h \ , \ k=0,1,\ldots,n \ .
\end{equation}
We will construct in each subinterval $[a+k h, a+(k+1)h]$ a
matrix-cubic spline approximating the solution of problem
(\ref{tema2}). For the first interval $[a,a+h] $, we consider that
the matrix-cubic spline is defined by

\begin{equation}\label{tema5}
S_{ \left|_{\left[a,a+h \right]}\right.}  (x) \ = \ Y(a) \ + \
Y'(a)(x-a) \ + \ \frac{1}{2!} Y''(a) (x-a)^2 \ + \ \frac{1}{3!}
A_{0} (x-a)^3 \ ,
\end{equation}

where the matrix $A_{0} \in {\mathbb C}^{r \times q}$ is a parameter
to be determined. It is straightforward to check:

$$
S_{ \left|_{\left[a,a+h \right]}\right.}  (a)=Y(a) \ , \quad
S'_{\left|_{\left[a,a+h \right]}\right.}  (a)=Y'(a)=f(a,Y(a)) \ .
$$

To fully determine the matrix-cubic spline we still must obtain
$Y''(a) $ and $A_{0}$.
 We consider the functions $h_{1}$ and $h_{2}$ defined by

$$
\begin{array}{c}
h_{1}:[a,b] \mapsto [a,b] \\
\\
h_{1}(x) \ = \ x
\end{array} \ , \qquad
\begin{array}{c}
h_{2}:[a,b] \mapsto {\mathbb C}^{r \times q} \\
\\
h_{2}(x) \ = \ Y(x)
\end{array} \ , \
$$

where $Y(x)$ is the theoretical solution of (\ref{tema2}). We
describe now $f(x, Y(x))$ as a composition of functions $f$ and
$\left(h_{1}, h_{2} \right)$, that is, let $\phi:[a,b] \
\mapsto \ {\mathbb C}^{r \times q}$ be defined by

$$
\phi(x) \ = \ \left[f \circ \left(h_{1}, h_{2}\right) \right](x) \ = \
f\left(h_{1}(x), h_{2}(x) \right) \ = \ f(x, Y(x)) \ .
$$

Thus, $\phi$ is a real variable function of $x$, and applying
theorem $8.9.2$ of \cite[p.170]{graham} its
derivative takes the form:

$$
D\phi \ = \ D\left(f \circ \left(h_{1}, h_{2}\right) \right) \ = \
\left(\left(D_{1} f \right) \left(h_{1}, h_{2}\right) \right)\cdot Dh_{1} \ + \
\left(\left(D_{2} f \right) \left(h_{1}, h_{2}\right) \right)\cdot Dh_{2} \ ,
$$

where the partial derivatives of $f$, $D_{1}(f)$, $D_{2}(f)$ exist
and are continuous since it is assumed that $f \in {\cal
C}^{1}\left(T \right)$. By (\ref{tema2}) it is clear that

$$
\frac{d \left(vec \ Y(x) \right)^{T}}{dx} \ = \
\left[vec \ f(x, Y(x))\right]^{T} \ .
$$

Next, applying the chain rule for matrix functions (\ref{tema0}) and
then taking the derivative of a matrix with respect to a matrix,
(\ref{tema1}), one obtains

\begin{equation}\label{tema6}
Y''(x) \ = \ \frac{\partial f(x, Y(x))}{\partial x} \ + \ \left[
\left[vec \ f(x, Y(x))\right]^{T} \otimes I_{r}\right]
\frac{\partial f(x, Y(x))}{\partial \ vec \ Y(x)} \ .
\end{equation}

We are now in the position to evaluate $Y''(a)$ using (\ref{tema6}). \\

By imposing that (\ref{tema5}) is a solution of problem (\ref{tema2}) in
$x=a+h$, we have:
\begin{equation}\label{3}
S'_{ \left|_{\left[a,a+h \right]}\right.} (a+h) \ = \ f\left(a+h,
S_{ \left|_{\left[a,a+h \right]}\right.}  (a+h) \right) \ ,
\end{equation}
and obtain from (\ref{3}) the matrix equation with only one unknown
matrix $A_{0}$:
\begin{equation}\label{4}
A_{0} \ = \ \frac{2}{h^2}\left[f\left(a+h,Y(a)+ Y'(a)h + \frac{1}{2} Y''(a) h^2 + \frac{1}{6} A_{0} h^3\right) - Y'(a) - Y''(a)
h \right] \ .
\end{equation}

Assuming that the matrix equation (\ref{4}) has only one solution
$A_{0}$, the
matrix-cubic spline is totally determined in the interval $[a,a+h]$. \\

Now, in the interval $[a+h, a+2h]$, the  matrix-cubic spline takes the
form

\begin{eqnarray}\label{tema10}
S_{ \left|_{\left[a+h, a+2h \right]}\right.}  (x) & = &
S_{ \left|_{\left[a, a+h \right]}\right.}  (a\!+\!h)\!+\!
S'_{ \left|_{\left[a, a+h \right]}\right.}  (a\!+\!h) (x-(a+h)) \nonumber \\
& + & \frac{1}{2!} S''_{ \left|_{\left[a, a+h \right]}\right.}  (a\!+\!h)
(x\!-\!(a\!+\!h))^2\!+\!\frac{1}{3!} A_{1} (x\!-\!(a\!+\!h))^3 \ ,
\end{eqnarray}
so that $S(x)$ is of class ${\cal C}^2([a,b]) $ on $[a,a+h] \cup
[a+h,a+2h]$, and all coefficients of the matrix-cubic spline $S_
{\left|_{\left[a+h, a+2h \right]} \right.} (x)$ are determined with
the exception of $A_{1} \in {\mathbb C}^{r \times q}$. By
construction, matrix-cubic spline (\ref{tema10}) satisfies the
differential equation (\ref{tema2}) in $x=a+h$. We can obtain
$A_{1}$ by requiring that the differential equation (\ref{tema2})
holds at point $x=a+2h$:

$$
S'_{ \left|_{\left[a+h,a+2h \right]}\right.} (a+2h) \ = \ f\left(a+2h,
S_{ \left|_{\left[a+h,a+2h \right]}\right.}  (a+2h) \right) \ .
$$

Expanding, we obtain the matrix equation with only one unknown
matrix $A_{1}$:

\begin{eqnarray}\label{4a}
A_{1} & = & \frac{2}{h^2}\left[f\left(a+2h,
S_{ \left|_{\left[a,a+h \right]}\right.}  (a+h)+
S'_{ \left|_{\left[a,a+h \right]}\right.}  (a+h) h +
\frac{1}{2} S''_{ \left|_{\left[a,a+h \right]}\right.}  (a+h) h^2 +
\frac{1}{6} A_{1} h^3\right) \right. \nonumber \\
& - & \left.
S'_{ \left|_{\left[a,a+h \right]}\right.}  (a+h) -
S''_{ \left|_{\left[a,a+h \right]}\right.}  (a+h)h \right] \ .
\end{eqnarray}

Let us assume that the matrix equation (\ref{4a}) has only one
solution $A_{1} $. This way the spline is totally determined
in the interval $[a+h,a+2h] $. \\

Iterating this process, let us construct the matrix-cubic spline
taking $\left[a+(k-1) h, a+k h \right]$ as the last subinterval. For
the next subinterval $\left[a+kh, a+(k+1)h\right]$, we define the
corresponding matrix-cubic spline as

\begin{equation}\label{tema14}
S_{ \left|_{\left[a+k h, a+(k+1)h \right]}\right.}  (x) =
\beta_{k}(x) \ + \ \frac{1}{3!} A_{k} (x-(a+k h))^3 \ ,
\end{equation}
where
\begin{eqnarray}\label{tema14c}
\beta_{k}(x) \ = \ \sum_{k=0}^{2} \frac{1}{k!}  S^{(k)}_{ \left|_{\left[a+(k-1)h,
a+kh \right]}\right.}  (a+k h) (x-(a+k h))^{k} \ .
\end{eqnarray}

With this definition, the matrix-cubic spline is $S(x) \in {\cal C}^2
\displaystyle \left (\bigcup_{j=0}^{k} [a+jh,a+(j+1)h] \right) $
and fulfills the differential equation
(\ref{tema2}) at point $x=a+kh $.
As an additional requirement, we assume that $S(x)$ satisfies the differential
equation (\ref{tema2}) at the point $x=a+(k+1)h$:

$$
S'_{ \left|_{\left[a+kh,a+(k+1)h \right]}\right.}
(a+(k+1)h) \ = \ f\left(a+(k+1)h, S_{ \left|_{
\left[a+kh,a+(k+1)h \right]}\right.}  (a+(k+1)h) \right) \ ,
$$

and expanding this equation with the unknown matrix $A_{k}$ yields
\begin{eqnarray}\label{tema15a}
A_{k} & = & \frac{2}{h^2} \left[ f\left( a+(k+1)h,
S_{ \left|_{\left[a+(k-1)h,a+kh \right]}\right.}  (a+kh)+
S'_{ \left|_{\left[a+(k-1)h,a+kh \right]}\right.}  (a+kh) h \right. \right. \nonumber \\
& + & \left. \frac{1}{2} S''_{ \left|_{\left[a+(k-1)h,a+kh \right]}\right.} (a+kh) h^2 +
\frac{1}{6} A_{k} h^3\right)  \nonumber \\
& - & \left.
S'_{ \left|_{\left[a+(k-1)h,a+kh \right]}\right.}  (a+kh) -
S''_{ \left|_{\left[a+(k-1)h,a+kh \right]}\right.}  (a+kh)h \right] \ .
\end{eqnarray}
Note that this matrix equation (\ref{tema15a}) is analogous to
equations (\ref{4}) and (\ref{4a}), when $k=0$ and $k=1$,
respectively. We will
show that these equations have an unique solution using a fixed-point argument.\\

For a fixed $h$, we will consider the matrix function of matrix
variable $g:{\mathbb C}^{r \times q} \mapsto {\mathbb C}^{r \times
q}$ defined by
\begin{eqnarray}\label{iterative}
g(T) & = & \frac{2}{h^2}\left[f\left(a+(k+1)h,
S_{ \left|_{\left[a+(k-1)h,a+kh \right]}\right.}  (a+kh)+
S'_{ \left|_{\left[a+(k-1)h,a+kh \right]}\right.}  (a+kh) h  \right.\right. \nonumber\\
& + & \left. \frac{1}{2} S''_{ \left|_{\left[a+(k-1)h,a+kh \right]}\right.}  (a+kh) h^2 +
\frac{1}{6} T h^3\right) \nonumber \\
& - & \left.
S'_{ \left|_{\left[a+(k-1)h,a+kh \right]}\right.}  (a+kh) -
S''_{ \left|_{\left[a+(k-1)h,a+kh \right]}\right.}  (a+kh)h \right] \ .
\end{eqnarray}
Relation (\ref{tema15a}) holds if and only if $A_{k}=g(A_{k})$,
that is, if $A_{k}$ is a
fixed point for function $g(T)$.\\

Observe that by using (\ref{tema14c}) and applying the global Lipschitz's
condition (\ref{2}) it follows that

$$
 \left\|g(T_{1}) - g(T_{2}) \right\| \leq \frac{Lh}{3} \left\|T_{1}-T_{2} \right\|
 .
$$

Taking $h<3/L$, $g(T)$ yields a contractive matrix function, which
guarantees that equation (\ref{tema15a}) has unique solutions
$A_{k}$ for $k=0,1,\ldots,n-1$.  Hence, the matrix-cubic spline is
completely determined. Taking into account \cite[Theorem 5]{20}, the
following result can be established.


\begin{theorem}\label{theorem1}
Let be $L$ the Lipschitz constant defined by (\ref{2}).  If $h \leq
3/L$, then the matrix-cubic spline $S(x)$ exists in each subinterval
$\left[a+kh, a+(k+1)h \right]$, $k=0,1,\ldots,n-1$, as defined in
the previous construction. Furthermore, if $f\in {\cal C}^{3}(T)$,
then $ \left\|Y(x) - S(x)\right\| = O(h^4) \ \forall x \in [a,b],$
where $Y(x) $ is the theoretical solution
of (\ref{tema2}).\\
\end{theorem}

\section{Algorithm}\label{seccion3}
The following algorithm is designed to compute the approximate
solution of (\ref{tema2}) by means of matrix-cubic splines in the
interval $[a,b]$ with an error of the order $O(h^4)$ under
conditions of theorem \ref{theorem1}.

\begin{center}
\framebox[13.5cm][l]{
\begin{tabular}{l}
$\bullet$ Determine the constant $Y ''(a)$ given by (\ref{tema6}).
Take $n>L(b-a)/3$, $h=(b-a)/n$ and \\
\hskip.3cm the partition $\Delta_{[a,\ b]}$ defined by Eq.~(\ref{partition}). \\[.3cm]

$\bullet$ Solve the matrix equation (\ref{4}) for $k=0$ and
determine $S_ {\left|_{\left[a,a+h \right]} \right.}(x)$
of Eq.~(\ref{tema5}). \\[.3cm]

$\bullet$ Solve the matrix equation (\ref{tema15a}) iteratively
for $k=1,\ldots,n-1 $, and then compute the \\
\hskip.3cm splines $S_ {\left|_{\left[a+k h,a+(k+1)h \right]} \right.}(x)$
according to Eq.~(\ref{tema14}).
\end{tabular}
   }
\end{center}

Depending on the function $f(t,Y)$, matrix equations (\ref{4}) and
(\ref{tema15a}) can be solved explicitly (see \cite{ult}) or using the
iterative method (see for example \cite{ort}):

$$
T^{s}_{l+1} \ = \ g(T^{s}_{l}) \ , \ \mbox{where} \ T^{s}_{0} \
\mbox{ is an arbitrary matrix in} \ {\mathbb C}^{r \times q} \ , \
s=0,1,\ldots,n-1
$$

and $g(T)$ is given for (\ref{iterative}).
In the following section, we will test the algorithm proposed.

\section{Example: A non-linear vector system}\label{seccion7}

We consider the next non-linear vector differential system

\begin{equation}\label{guzman1}
\left.  \begin{array}{rcl}
y '_{1}(x) & = & -1 + e^x - \sin{(x)} + \sin{(y_2(x))}  \\
\\
y '_{2}(x) & = & \frac{1}{4 + y_1(x)^2} - \frac{1}{5 + e^{2\,x} +
2\,e^x\,\cos{(x)} - \sin^2{(x)}} \\
\\
y_{1}(0) & = & 2 , \qquad y_{2}(0) \ = \ \frac{\pi}{2}
\end{array}
 \right\}  \ 0 \leq x \leq 1 \ ,
\end{equation}

It is easy to check that this problem has the exact solution
$y_{1}(x)=e^{x}+\cos{(x)}, y_{2}(x)=\frac{\pi}{2}$, so in this
particular case we will be able
to obtain the exact error of our numerical estimates.\\

We can rewrite (\ref{guzman1}) in the compact form
\begin{equation}\label{guzman2}
\left.  \begin{array}{rcl}
Y'(x) & = & F\left(x,Y\right)  \\
\\
Y(0)  & = & \left(\begin{array}{c} 2\\  \frac{\pi}{2}\end{array}
\right)
\end{array}
 \right\}  \ 0 \leq x \leq 1 \ , \ Y(x)=\left(\begin{array}{c} y_{1}(x)\\
y_{2}(x)\end{array} \right) \in {\mathbb R}^{2} \ , \  F(x,Y)=\left(\begin{array}{c} -1 + e^x - \sin{(x)} + \sin{(y_2(x))}\\
 \frac{1}{4 + y_1(x)^2} - \frac{1}{5 + e^{2\,x} +
2\,e^x\,\cos{(x)} - \sin^2{(x)}}\end{array} \right) \in {\mathbb
R}^{2},
\end{equation}

thus $Y'(0)=F\left(0,\left(\begin{array}{c} 2\\
\frac{\pi}{2}\end{array} \right)\right)=\left(\begin{array}{c} 1\\
0\end{array} \right)$. We calculate $Y''(0)$ using (\ref{tema6}). in
this case, one gets

\begin{equation}\label{guzman12}
vec(Y(x))=Y(x)=\left(\begin{array}{c} y_{1}(x)\\
y_{2}(x)\end{array} \right) \ , \ \frac{\partial F(x,
Y(x))}{\partial x} = \left(\begin{array}{c}
 e^x - \cos{(x)} \\
 \\
\frac{2 e^{2 x} + 2 e^x \cos{(x)} -
     2 e^x \sin{(x)} - 2\cos{(x)}\sin{(x)}}{
     \left( 5 + e^{2 x} +
        2 e^x \cos{(x)} - \sin^2{(x)} \right)^2
        }
        \end{array} \right) .
\end{equation}
On the other hand, we have
$$
\left[vec \ F(x, Y(x))\right]^{T} \otimes I_{2}
$$
\begin{eqnarray*}
&=&\left(\begin{array}{cc} -1 + e^x - \sin{(x)} + \sin{(y_2(x))} &
 \frac{1}{4 + y_1(x)^2} - \frac{1}{5 + e^{2x} +
2e^x\cos{(x)} - \sin^2{(x)}}\end{array} \right) \otimes I_{2} \\
&=& \left(\begin{array}{c|c} \left(-1 + e^x - \sin{(x)} +
\sin{(y_2(x))}\right)I_2 &
 \left(\frac{1}{4 + y_1(x)^2} - \frac{1}{5 + e^{2x} +
2e^x\cos{(x)} - \sin^2{(x)}}\right)I_2\end{array} \right)\nonumber\\
\end{eqnarray*}
\begin{equation}\label{guzman13}
=\!\!{\footnotesize \left(\begin{array}{cc|cc}
 -\!1\!+\!e^x\!-\!\sin{(x)}\!+\!\sin{(y_2(x))} & 0 & \frac{1}{4 + y_1(x)^2} - \frac{1}{5 + e^{2x}
+ 2e^x\cos{(x)} - \sin^2{(x)}} & 0\\
0 &  -\!1\!+\!e^x\!-\!\sin{(x)}\!+\!\sin{(y_2(x))} & 0 & \frac{1}{4
+ y_1(x)^2} - \frac{1}{5 + e^{2x} + 2e^x\cos{(x)} - \sin^2{(x)}}
 \end{array} \right)},
\end{equation}
and
\begin{equation}\label{guzman14}
\frac{\partial F(x, Y(x))}{\partial \ vec \ Y(x)}=
\left(\begin{array}{c}
\frac{\partial F(x, Y(x))}{\partial y_1 } \\
\\
\hline \\
\frac{\partial F(x, Y(x))}{\partial y_2}
\end{array}
\right)=
\left(\begin{array}{c} \frac{\partial}{\partial y_1 } \left( -1 + e^x - \sin{(x)} + \sin{(y_2(x))} \right) \\
\frac{\partial}{\partial y_1 }\left( \frac{1}{4 + y_1(x)^2} -
\frac{1}{5 + e^{2\,x} + 2\,e^x\,\cos{(x)} - \sin^2{(x)}} \right) \\
\hline \\
\frac{\partial}{\partial y_2 } \left( -1 + e^x - \sin{(x)} + \sin{(y_2(x))} \right) \\
\frac{\partial}{\partial y_2 }\left( \frac{1}{4 + y_1(x)^2} -
\frac{1}{5 + e^{2\,x} + 2\,e^x\,\cos{(x)} - \sin^2{(x)}} \right) \\
 \end{array}
\right)=
\left(\begin{array}{c} 0 \\
\frac{-2 y_1(x)}{\left( 4 + y_1(x)^2 \right)^2} \\
\hline \\
\cos{(y_2(x))} \\
0 \\
 \end{array}
\right) ,
\end{equation}

Therefore

\begin{equation}\label{guzman15}
\left[ \left[vec \ F(x, Y(x))\right]^{T} \otimes
I_{2}\right]\frac{\partial F(x, Y(x))}{\partial \ vec \ Y(x)}=
\left(\begin{array}{c} \left( \frac{1}{4 + y_1(x)^2} -
     \frac{2}{9 + 2e^{2x} + 4e^x\cos{(x)} + \cos{(2x)}} \right) \,
   \cos{(y_2(x))} \\
   -\frac{2 y_1(x)\left( -1 + e^x - \sin{(x)} + \sin{(y_2(x))} \right)}
   {\left( 4 +y_1(x)^2 \right)^2} \end{array} \right),
\end{equation}

and by (\ref{guzman12})-(\ref{guzman15}) one concludes

\begin{eqnarray}\label{guzman4}
Y''(x) & = & \frac{\partial F(x, Y(x))}{\partial x} \ + \ \left[
\left[vec \ F(x, Y(x))\right]^{T} \otimes I_{2}\right]
\frac{\partial F(x, Y(x))}{\partial \ vec \ Y(x)} \nonumber \\
& = & \left(\begin{array}{c}  e^x - \cos{(x)}+ \left( \frac{1}{4 +
y_1(x)^2} -
     \frac{2}{9 + 2e^{2x} + 4e^x\cos{(x)} + \cos{(2x)}} \right) \,
   \cos{(y_2(x))} \\
\frac{2 e^{2 x} + 2 e^x \cos{(x)} -
     2 e^x \sin{(x)} - 2\cos{(x)}\sin{(x)}}{
     \left( 5 + e^{2 x} +
        2 e^x \cos{(x)} - \sin^2{(x)} \right)^2
        }-\frac{2 y_1(x)\left( -1 + e^x - \sin{(x)} + \sin{(y_2(x))} \right)}
   {\left( 4 +y_1(x)^2 \right)^2} \end{array} \right).
\end{eqnarray}

Taking into account that $y_1(0)=2,y_{2}(0)=\frac{\pi}{2}$ and
evaluating $Y''(x)$ of (\ref{guzman4}) when $x=0$, one gets
$Y''(0)=\left(\begin{array}{c}  0 \\
0 \end{array} \right)$.\\

It is straightforward to show that $F$, defined by (\ref{guzman2}),
fulfills the global Lipschitz's condition
\begin{equation}\label{guzman3}
\left\|f\left(x, Y\right) \ - \ f\left(x, Z\right)\right\| \ \leq
\left\|Y - Z \right\| \ , \  0 \leq x \leq 1 \ , Y, Z \in {\mathbb
R}^{2} \ ,
\end{equation}
thus, we can take $L$ given by (\ref{2}) as $L=1$. Therefore, we
need to take $h<3/L$ and thus $h=0.1$ for example. The results are
generated with {\sl Mathematica} using {\tt FindRoot} function to
solve the emerging algebraic equations, and are summarized in Table
1. In each interval, we evaluated the difference between the
estimates of our numerical approach and the exact solution, and then
take the Frobenius norm of this difference. The maximum of these
errors are indicated in the third column for
each subinterval.\\

\begin{center}
\begin{table}
{\small
\begin{tabular}{|c|c|c|}
\hline
\mbox{Interval} & \mbox{Approximation} & \mbox{Max.\ Error} \\
\hline $[0,0.1]$ & {\small $\left( \begin{array}{c}
2 + x + 0.177917 x^3  \\
\frac{\pi }{2} - 5.62424 \times 10^{-6} x^3
\end{array}
\right)$ } & {\small $2.83337 \times 10^{-6}$} \\
\hline $[0.1,0.2]$ & {\small $\left(
\begin{array}{c}
1.99995 + 1.00138 x - 0.0138342 x^2 + 0.224031 x^3 \\
1.5708 + 6.67857 \times 10^{-7} x - 6.67857\times 10^{-6} x^2 +
    0.0000166377 x^3 \end{array}
\right)$ } & {\small $ 2.83337 \times 10^{-6}$}\\
\hline $[0.2,0.3]$ & {\small $\left(
\begin{array}{c}
1.99975 + 1.00445 x - 0.0291822 x^2 + 0.249611 x^3\\
1.5708 - 4.57386\times 10^{-6}x + 0.00001953x^2 - 0.0000270433x^3
\end{array}
\right)$ } & {\small $2.94712 \times 10^{-6}$ } \\
\hline $[0.3,0.4]$ & {\small $\left(
\begin{array}{c}
1.99841 + 1.01783 x - 0.0737602 x^2 + 0.299142 x^3 \\
1.57079 + 0.0000126873 x - 0.0000380073 x^2 + 0.0000368871 x^3
\end{array}
\right)$ } & {\small $ 2.94712 \times 10^{-6}$ } \\
\hline $[0.4,0.5]$ & {\small $\left(
\begin{array}{c}
1.99655 + 1.0318 x - 0.108685 x^2 + 0.328246 x^3 \\
1.5708 - 0.0000271633 x + 0.0000616192 x^2 - 0.0000461351 x^3
\end{array}
\right)$ } & {\small $3.0698 \times 10^{-6}$} \\
\hline $[0.5,0.6]$ & {\small $\left(
\begin{array}{c}
1.9899 + 1.07171 x - 0.188509 x^2 + 0.381462 x^3 \\
1.57079 + 0.0000485499 x - 0.0000898071 x^2 + 0.0000548159 x^3
\end{array}
\right)$ } & {\small $3.0698 \times 10^{-6}$} \\
\hline $[0.6,0.7]$ & {\small $\left(
\begin{array}{c}
1.98277 + 1.10736 x - 0.247933 x^2 + 0.414475 x^3 \\
1.57081 - 0.0000785694 x + 0.000122058 x^2 - 0.0000628872 x^3
\end{array}
\right)$ } & {\small $3.20977 \times 10^{-6}$ } \\
\hline $[0.7,0.8]$ & {\small $\left(
\begin{array}{c}
1.96307 + 1.19177 x - 0.368517 x^2 + 0.471896 x^3 \\
1.57077 + 0.000117333 x - 0.000157802 x^2 + 0.0000703796 x^3
\end{array}
\right)$ } & {\small $3.20977 \times 10^{-6}$ } \\
\hline $[0.8,0.9]$ & {\small $\left(
\begin{array}{c}
1.94382 + 1.26395 x - 0.45874 x^2 + 0.509489 x^3\\
1.57084 - 0.0001661 x + 0.00019649 x^2 - 0.0000772419 x^3
\end{array}
\right)$ } & {\small $3.37764 \times 10^{-6}$ } \\
\hline $[0.9,1.0]$ & {\small $\left(
\begin{array}{c}
1.89829 + 1.41574 x - 0.627395 x^2 + 0.571954 x^3 \\
1.57073 + 0.000224533 x - 0.000237548 x^2 + 0.0000835127 x^3
\end{array}
\right)$ } & {\small $3.37764 \times 10^{-6}$ } \\
\hline
\end{tabular}
} \caption{Approximation for vector differential
system~(\ref{guzman1}) in the interval $[0,1]$ with step size
$h=0.1$.}
\end{table}
\end{center}

\begin{center}
\begin{figure}
\includegraphics[scale=1.3]{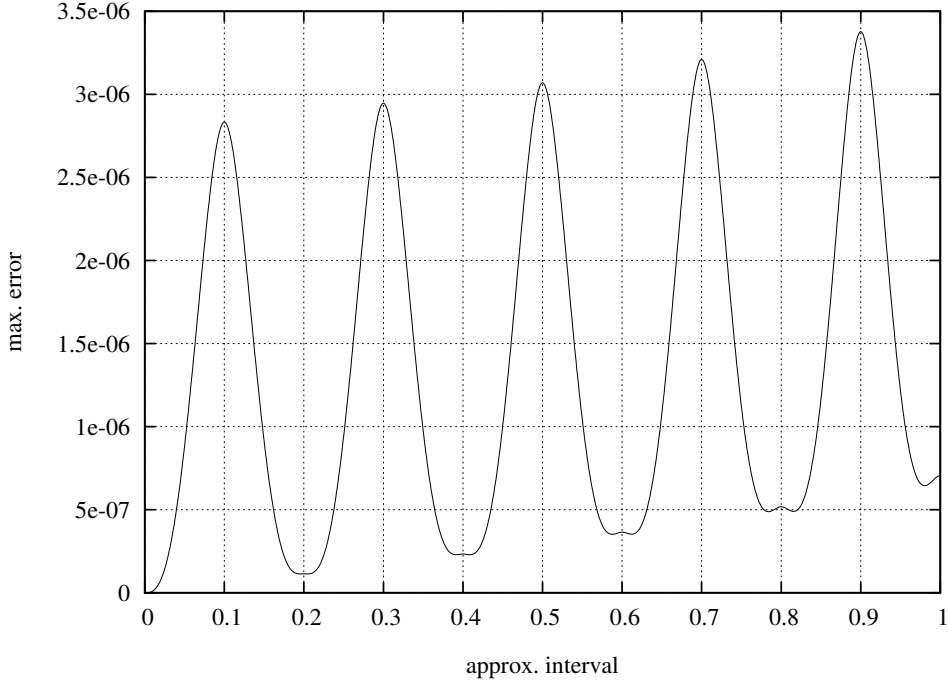}
\caption{Representing the Frobenius error margins for vector
differential system~(\ref{guzman1}) in the interval $[0,1]$ with
step size $h=0.1$.}
\end{figure}
\end{center}

\section{Example: Sylvester matrix differential equation}\label{seccion4}

Linear matrix differential equations of the type

\begin{equation}\label{sylvester1}
\left.  \begin{array}{rcl}
Y'(x) & = & A(x)Y(x)+Y(x)B(x)+C(x)  \\
\\
Y(a)  & = & Y_{a}
\end{array}
 \right\}  \ a \leq x \leq b \ , \ Y(x),A(x), B(x), C(x) \in {\mathbb C}^{r \times
 r} \ ,
\end{equation}

arise in many fields of science and engineering. In the case of
constant coefficients has been studied by several authors (see for
example \cite{locas1}). However, the variable-coefficient case
has so far received little
numerical treatment in the literature. We can observe that the
proposed method require the matrix functions $A(x), B(x)$ and $C(x)$
to be differentiable, while, for example, in the method proposed in
\cite{locas}, it is necessary that $A(x), B(x)$ have continuous
second-order derivatives and $C(x)$ continuous in the domain $a \leq
x \leq b$. \\

As an example, here let us consider the Sylvester problem~(\ref{sylvester1})
with
\begin{equation}\label{paramsylv}
\left.
\begin{array}{lcl}
   A(x) & = & \left(
   \begin{array}{cc}
          0 & xe^{-x} \\
          x & 0
    \end{array}\right), \quad
B(x) =
\left(\begin{array}{cc}
0 & x \\
0 & 0
\end{array} \right), \quad
C(x) =
\left(\begin{array}{cc}
-e^{-x}(1+x^2) & -2e^{-x}x \\
1-e^{-x}x & -x^2
 \end{array} \right) \\
 \\
Y(0) &=&
\left(\begin{array}{rr}
1 & 0\\
0 & 1
\end{array}
\right), \quad
Y(x) \in {\mathbb C}^{2 \times 2}  \ , \ 0 \leq x \leq 1
\end{array}  \right\}
\end{equation}

This problem has an exact solution $ Y(x)=\left(
\begin{array}{cc}
e^{-x} & 0  \\
x & 1
\end{array}
\right)$, so in this particular case we will be able to obtain the
exact error of our numerical estimates.\\

As we have $\displaystyle\max_{x\in [0,1]}\left(\left\|
\left(\begin{array}{cc}
0 & xe^{-x} \\
x & 0
\end{array}
\right)\right\|+\left\| \left(\begin{array}{cc}
0 & x \\
0 & 0
\end{array} \right)\right\|\right)
 \leq 1.69443$, one can take the constant $L$ given for (\ref{2}) as $L=2$.\\


Taking the derivative of $Y'(x)=A(x)Y(x)+Y(x)B(x)+C(x)$, gives:

\begin{eqnarray}\label{tema6a}
Y''(x) & = &
\left(A'(x) + \left(A(x)\right)^2\right)Y(x) + Y(x)\left(\left(B(x)\right)^2+B'(x)\right)\nonumber \\
\quad &+& 2A(x)Y(x)B(x)+ A(x)C(x)+C(x)B(x)+C'(x) \ .
\end{eqnarray}
We see that $Y'(0) \ = \ \left(
\begin{array}{cc}
-1 & 0  \\
1 & 0
\end{array}
\right)$, and by applying (\ref{tema6a}) it is $Y''(0) \ = \ \left(
\begin{array}{cc}
1 & 0  \\
0 & 0
\end{array}
\right)$.\\

In this numerical example, we take $n=10 $ such that $n>L(b-a)/3$
and $h=0.1=(b-a)/n$. The results are generated with {\sl
Mathematica} using the Bartels-Stewart algorithm (see for example
\cite{4}) to solve the emerging algebraic equations, and are
summarized in Table 2, where the numerical estimates have been
rounded to the fourth relevant digit. In each interval, we evaluated
the difference between the estimates of our numerical approach and
the exact solution, and then take the Frobenius norm of this
difference. The maximum of these errors are indicated in the third
column for each subinterval.

\begin{center}
\begin{table}
{\small
\begin{tabular}{|c|c|c|}
\hline
\mbox{Interval} & \mbox{Approximation} & \mbox{Max.\ Error} \\
\hline $[0,0.1]$ & {\small $\left( \begin{array}{cc}
1 - x + 0.5 x^2 - 0.1612 x^3 & 0 \\
x & 1
\end{array}
\right)$ } & {\small $1.33472 \times 10^{-6}$} \\
\hline $[0.1,0.2]$ & {\small $\left(
\begin{array}{cc}
1 - 0.9994 x + 0.4938 x^2 - 0.1406 x^3 & 0 \\
x & 1 \end{array}
\right)$ } & {\small $1.33472 \times 10^{-6}$}\\
\hline $[0.2,0.3]$ & {\small $\left(
\begin{array}{cc}
1 - 0.9984 x + 0.4890 x^2 - 0.1325 x^3 & 0 \\
x  & 1 \end{array}
\right)$ } & {\small $1.2445 \times 10^{-6}$ } \\
\hline $[0.3,0.4]$ & {\small $\left(
\begin{array}{cc}
0.9994 - 0.9936 x + 0.4728 x^2 - 0.1146 x^3 & 0 \\
x  & 1
\end{array}
\right)$ } & {\small $1.2445 \times 10^{-6}$ } \\
\hline $[0.4,0.5]$ & {\small $\left(
\begin{array}{cc}
0.9991-0.9909x+0.4661x^{2}-0.1090x^{3} & -0.0001 x^2 \\
1.0001x-0.0001x^2 & 1
\end{array}
\right)$ } & {\small $1.17402 \times 10^{-6}$} \\
\hline $[0.5,0.6]$ & {\small $\left(
\begin{array}{cc}
0.9971-0.9791x+0.4426x^{2}-0.0933x^{3} & -0.0001x+0.0002x^2-0.0001x^3 \\
0.9999x+0.0002x^2 - 0.0001x^3 & 0.9999
\end{array}
\right)$ } & {\small $1.17402 \times 10^{-6}$} \\
\hline $[0.6,0.7]$ & {\small $\left(
\begin{array}{cc}
0.9963-0.9732x+0.4361x^{2}-0.0898x^{3} & 0.0002x-0.0004x^2 +0.0002x^3  \\
1.0002x-0.0004x^2 + 0.0002x^3 & 1
\end{array}
\right)$ } & {\small $1.12331 \times 10^{-6}$ } \\
\hline $[0.7,0.8]$ & {\small $\left(
\begin{array}{cc}
0.9916-0.9549x+0.4071x^{2}-0.07591x^{3} & 0.0001-0.0004x+0.0006x^2 - 0.0003x^3 \\
0.0001+0.9996x+0.0006x^2 -0.0003x^3 & 0.9999
\end{array}
\right)$ } & {\small $1.12331 \times 10^{-6}$ } \\
\hline $[0.8,0.9]$ & {\small $\left(
\begin{array}{cc}
0.9906-0.9512x+0.4025x^{2}-0.0739x^{3} & 0.0002+0.0007x-0.0009x^2 +0.0003x^3 \\
-0.0002+1.0007x-0.0009x^2 +0.0003x^3 & 1+0.0001x^2
\end{array}
\right)$ } & {\small $1.09412 \times 10^{-6}$ } \\
\hline $[0.9,1.0]$ & {\small $\left(
\begin{array}{cc}
0.9816-0.9212x+0.3691x^{2}-0.0616x^{3} & 0.0004-0.0011x+0.0012x^2 -0.0004x^3 \\
0.0004+0.0.9989x+0.0012x^2 - 0.0004x^3 & 0.9999+0.0002x-0.0002x^2
\end{array}
\right)$ } & {\small $1.09412 \times 10^{-6}$ } \\
\hline
\end{tabular}
} \caption{Approximation for Sylvester matrix differential
equation~(\ref{paramsylv}) in the interval $[0,1]$ with step size
$h=0.1$.}
\end{table}
\end{center}

\begin{center}
\begin{figure}
\includegraphics[scale=1.3]{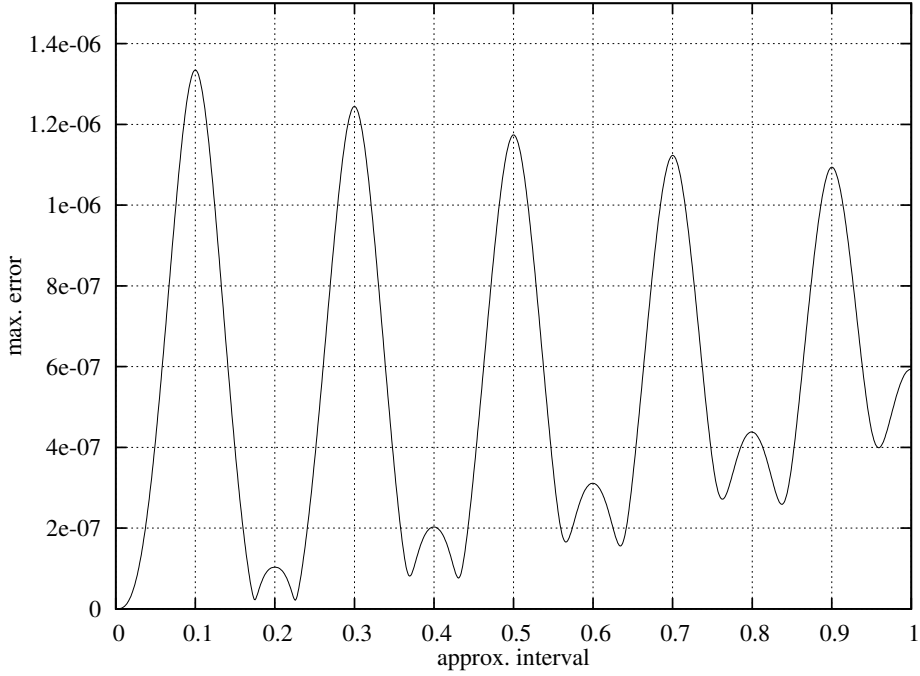}
\caption{Representing the absolute error margins for the Sylvester
matrix differential equation~(\ref{paramsylv}) in the interval $[0,1]$
with step size $h=0.1$.}
\end{figure}
\end{center}

\section{Example: Riccati matrix differential equation}\label{seccion5}

Rectangular non-symmetric Riccati matrix-differential
equation of the type
\begin{equation}\label{riccati1}
\left.  \begin{array}{rcl}
Y'(x) & = & C(x)-D(x)Y(x)-Y(x)A(x)-Y(x)B(x)Y(x)  \\
\\
Y(0)  & = & Y_{0}
\end{array}
 \right\} \quad 0 \leq x \leq c,
\end{equation}
where the unknown $Y(x) \in {\Bbb C}^{p \times q}$ and coefficients
$A(t) \in {\Bbb C}^{q \times q}, B(x) \in {\Bbb C}^{q \times p},
C(x) \in {\Bbb C}^{p \times q}, D(x) \in {\Bbb C}^{p \times p}$ are
differentiable matrix-valued functions
arise frequently in important
applications to classical control theory \cite{casti} and as
decoupling techniques for both the analytic and numerical study of
boundary value problems \cite{Ascher}.
The Riccati equation (\ref{riccati1}) has been studied extensively,
and different resolution techniques have been introduced (see
\cite{Jodar3} and references therein).\\

The study of the Riccati equation (\ref{riccati1}) is closely
related to the underlying linear system
\begin{equation}\label{riccati2}
\left.  \begin{array}{rcl}
X'(x) & = & S(x) X(x) \\
\\
X(0)  & = & \left[\begin{array}{c} I_{q} \\
Y_{0} \end{array} \right]
\end{array}
 \right\}  \ \mbox{where} \ X(x) \ = \ \left[\begin{array}{c}
U(x) \\
V(x)
  \end{array} \right] \ , \ S(x) \ = \ \left[\begin{array}{rr}
A(x) & B(x) \\
C(x) & -D(x) \end{array} \right]\ .
\end{equation}
Specifying the solution of (\ref{riccati1}) is given by
\begin{equation}\label{riccati3}
Y(x) \ = \ V(x) U^{-1}(x)
\end{equation}
where $Y(x)$ is defined in the interval where $U(x)$ is
invertible, see \cite{5}.\\

Taking into account lemma 1 and 2 of \cite{Jodar2}, $U(x)$ is
invertible in the interval $[0, \delta]$ and the solution $Y(x)$ of
problem (\ref{riccati1}) satisfies
\begin{equation}\label{ricca4}
\left\|Y(x)\right\| \ \leq \ M \ , \ M \ = \ \left(1-\delta q_{0}
\exp{(\delta k_{0})} w_{0} \right)^{-1} w_{0} \exp{(\delta k_{0})},
\end{equation}
where $\delta$ is a positive number satisfying
\begin{equation}\label{ricca5a}
\delta k_{0} + \log{(\delta)} < -\log{(q_{0} w_{0})} \ ,
\end{equation}
and
\begin{equation}\label{ricca5}
\left.\begin{array}{rcl} k_{0} & = &
\max\left\{\left\|\left[\begin{array}{rr}
A(x) & B(x) \\
C(x) & -D(x) \end{array} \right] \right\| \ ; \ 0 \leq x \leq
c\right\}  \\[.5cm]
q_{0} & = & \max\left\{\left\|A(x) \ B(x) \right\| \ ; \
0 \leq x \leq c\right\} \\[.5cm]
w_{0} & = & \left\|\begin{array}{c} I_{q} \\
Y_{0} \end{array} \right\|
\end{array} \right\}
\end{equation}

In accordance with \cite[p.1064]{locas2}, we consider the matrix-valued
function
\begin{equation}\label{ricca2}
F(x,Y) \ = \ C(x)-D(x) Y - Y A(x)-Y B(x) Y \ ,
\end{equation}
then, if we define
\begin{equation}\label{ricca1}
\left.\begin{array}{rcl} a & = & \sup\left\{\left\|A(x)\right\| ;
0\leq x \leq \delta \right\} \\
b & = & \sup\left\{\left\|B(x)\right\| ; 0\leq x \leq \delta
\right\} \\
c & = & \sup\left\{\left\|C(x)\right\| ; 0\leq x \leq \delta
\right\} \\
d & = & \sup\left\{\left\|D(x)\right\| ; 0\leq x \leq \delta
\right\}
\end{array}\right\}
\end{equation}
and $\left\|Y\right\| \leq M, \left\|\widetilde{Y}\right\| \leq M$,
with $M$ gives by (\ref{ricca4}), the following local Lipschitz
condition holds
\begin{equation}\label{local2}
\left\|F(x, Y) - F(x, \widetilde{Y}) \right\| \leq L
\left\|Y-\widetilde{Y} \right\| \ , \ L \ = \ a+d+2b M \ .
\end{equation}

In addition , if $\left\|Y\right\| \leq N$,
\begin{equation}\label{nueva}
\left\|F(x,Y) \right\| \leq c+N(a+d+bN) \ .
\end{equation}

Using the proposed spline method, the only one solution of
the matrix equations (\ref{4}) and (\ref{tema15a}) for
$k=1,\ldots,n-1$ is guaranteed using a fixed-point argument and
the global Lipschitz's condition (\ref{2}). In our case, we need to
prove the only one solution of the matrix equations (\ref{4}) and
(\ref{tema15a}) using a fixed point argument
and the local Lipschitz's condition (\ref{local2}).\\

We start with the matrix equation~(\ref{4}).
Let us suppose that $\left\|T \right\| \leq  N_1$. Taking into
account (\ref{nueva}), we take
\begin{equation}\label{nue5}
\left\{
\begin{array}{rcl}
N_2 & = & \left\|Y(a)\right\|+h\left\|Y'(a)\right\|+ \frac{h^2}{2}
\left\|Y''(a)\right\|+\frac{h^3}{6}N_1 \\
\ \ \\
N_3 & = & c+N_2(a+d+bN_2) \\
\ \ \\
N_4 & = & \frac{2}{h^2}\left(N_3+ \left\|Y'(a)\right\|+h
\left\|Y''(a)\right\|\right)
\end{array}  \right.
\end{equation}
with $a, b, c$ given by (\ref{ricca1}), and let be
$N=\max\left\{N_1,N_2,N_3,N_4,M \right\}$ with $M$ gives by
(\ref{ricca4}). Let be ${\cal A}=\left\{Y \in {\Bbb C}^{r \times q}
; \left\|Y\right\|\leq N \right\}$ and we consider the continuous
matrix-valued function of matrix variable $g:{\Bbb C}^{r \times q}
\mapsto {\Bbb C}^{r \times q}$ defined by~(\ref{iterative}) for $k=0$.
It is simple to verify that if $T \in {\cal A}$, by (\ref{nue5}) and
(\ref{nueva}) then $g(T) \in {\cal A}$. Thus, $g:{\cal A} \mapsto
{\cal A}$ and $A_{0}$ is a fixed point of $g$. In addition, if $T_1,
T_2 \in {\cal A} \ , \ \left\|T_1\right\|\leq M,
\left\|T_2\right\|\leq M$, has then that for $f$ defined by
(\ref{ricca2}), $f$ fulfills the local Lipschitz's condition (\ref{local2}) and
taking $h<3/L$, $g(T)$ yields a contractive matrix function, which
guarantees that equation (\ref{4}) has unique solutions $A_{0}$.
Hence, the matrix-cubic spline is completely determined in $[a,a+h]$.\\

For $k=1,\ldots,n-1$, fixed, supposed construct cubic-matrix spline
$S(x)$ taking $\left[a+(k-1) h, a+k h \right]$ as the last
subinterval, for the next subinterval $\left[a+kh, a+(k+1)h\right]$,
to define the corresponding spline we need determine $A_{k} \in
{\Bbb C}^{r \times q}$ as the only one solution of the matrix
equation~(\ref{tema15a}).
Let us suppose that $\left\|T \right\| \leq  \widetilde{N_1}$.
Taking into account (\ref{nueva}), we take
\begin{equation}\label{nue5y}
\left\{
\begin{array}{rcl}
\widetilde{N_2} & = & \left\|S_{ \left|_{\left[a+(k-1)h,a+kh
\right]}\right.} (a\!+\!kh)\right\|\!+\!h\left\|S'_{
\left|_{\left[a+(k-1)h,a+kh \right]}\right.} (a\!+\!kh)\right\|\!+\!
\frac{h^2}{2} \left\|S''_{ \left|_{\left[a+(k-1)h,a+kh
\right]}\right.} (a\!+\!kh)\right\|\!+\!\frac{h^3}{6}\widetilde{N_1} \\
\ \ \\
\widetilde{N_3} & = & c+\widetilde{N_2}(a+d+b\widetilde{N_2}) \\
\ \ \\
\widetilde{N_4} & = & \frac{2}{h^2}\left(\widetilde{N_3}+\left\|S'_{
\left|_{\left[a+(k-1)h,a+kh \right]}\right.}  (a\!+\!kh)\right\|+h
\left\|S''_{ \left|_{\left[a+(k-1)h,a+kh \right]}\right.}
(a\!+\!kh)\right\|\right)
\end{array}  \right.
\end{equation}
with $a, b, c$ given by (\ref{ricca1}), and let be
$\widetilde{N}=\max\left\{\widetilde{N_1},\widetilde{N_2},\widetilde{N_3},\widetilde{N_4},M
\right\}$ with $M$ gives by (\ref{ricca4}). Let be ${\cal
A}=\left\{Y \in {\Bbb C}^{r \times q} ; \left\|Y\right\|\leq
\widetilde{N} \right\}$ and we consider the continuous matrix-valued
function of matrix variable $g:{\Bbb C}^{r \times q} \mapsto {\Bbb
C}^{r \times q}$ defined by~(\ref{iterative}).

It is simple to verify that if $T \in {\cal A}$, by (\ref{nue5y})
and (\ref{nueva}) then $g(T) \in {\cal A}$. Thus, $g:{\cal A}
\mapsto {\cal A}$ and $A_{k}$ is a fixed point of $g$. In addition,
if $T_1, T_2 \in {\cal A} \ , \ \left\|T_1\right\|\leq M,
\left\|T_2\right\|\leq M$, has then that for $f$ defined by
(\ref{ricca2}), $f$ fulfills the local Lipschitz's condition (\ref{local2}) and
taking $h<3/L$, $g(T)$ yields a contractive matrix function, which
guarantees that equation (\ref{tema15a}) has unique solutions
$A_{k}$. Hence, the matrix-cubic spline is completely determined.\\

As an additional example for our proposed method, we consider the
Riccati matrix differential equation~(\ref{riccati1}) with
\begin{eqnarray}\label{paramricc}
A(x) &=& \left(
\begin{tabular}{ll}
$-x$ & $0$ \\
$-x$ & $x$
\end{tabular}
\right),\
B(x) = \left(
\begin{tabular}{ll}
$-x^2$ & $-2$ \\
$0$ & $1$
\end{tabular}
\right),\ D(x) \ = \ \left(
\begin{tabular}{ll}
$-1$ & $-x^{2}$ \\
$x$ & $x$
\end{tabular}
\right), \nonumber\\
C(x) & = & \left(
\begin{tabular}{ll}
$x\left( -e^{x}+e^{x}x-x^{3}\right) $ & $x\left( 2e^{x}-x^{2}\right) $ \\
$(1-x)x(2+x+2x^2)$ & $1+\left(
3-2x\right)
x^{2}+e^{x}\left( x-x^{4}\right)$
\end{tabular}
\right),\
 Y\left( 0\right) =\left(
\begin{tabular}{ll}
$0$ & $1$ \\
$0$ & $0$
\end{tabular}
\right).
\end{eqnarray}

In this case, the problem has an exact solution given by
$Y(x)=\left(
\begin{tabular}{ll}
$0$ & $e^{x}$ \\
$x^{2}$ & $x$
\end{tabular}
\right),$
which will permit us to obtain the total error for all our
numerical estimates.
A short computation using expressions~(\ref{ricca4})--(\ref{local2})
yields the following constants
\begin{equation}
\left.
\begin{array}{ll}
k_{0} = 6.13866 & \qquad
q_{0} = 3 \\
w_{0} = \sqrt{2} & \qquad
\delta = 0.115758 \\
M = 12.0883 & \qquad
a = 0.173205 \\
b = 2.23609 & \qquad
c = 1.17928 \\
d = 1.01 & \qquad
L = 55.2443
\end{array}
\right\}
\end{equation}
which are necessary for the spline approximation in the interval
$[0,0.1]$, where $\delta=0.1$ is taken for convenience. Therefore,
we need to take $h<3/L=0.0543042$ and thus $h=0.01$. The results are
generated with {\sl Mathematica} using {\tt FindRoot} function to
solve the emerging algebraic equations, and are summarized in Table
3, where the numerical estimates have been rounded to the fourth
relevant digit. In each interval, we evaluated the difference
between the estimates of our numerical approach and the exact
solution, and then take the Frobenius norm of this difference. The
maximum of these errors are indicated in the third column for each
subinterval.


\begin{center}
\begin{table}
{\small
\begin{tabular}{|c|c|c|}
\hline
\mbox{Interval} & \mbox{Approximation} & \mbox{Max.\ Error} \\
\hline
$[0,0.01]$ &
$\left(\begin{array}{cc}
0 & 1 + x + 0.5\,x^2 + 0.167224\,x^3 \\ x^2 & x
\end{array}\right)$
& $1.39903\times 10^{-10}$ \\
\hline
$[0.01,0.02]$ &
$\left(\begin{array}{cc}
0 & 1. + x + 0.499933\,x^2 + 0.169461\,x^3 \\ x^2 & x
\end{array}\right)$
& $1.39903\times 10^{-10}$ \\
\hline
$[0.02,0.03]$ &
$\left(\begin{array}{cc}
0 & 1 + x + 0.499864\,x^2 + 0.17061\,x^3 \\ x^2 & x
   \end{array}\right)$
& $ 1.41977\times 10^{-10}$ \\
\hline
$[0.03,0.04]$ &
$\left(\begin{array}{cc}
0 & 1 + 1.00001\,x + 0.49966\,x^2 + 0.172877\,x^3 \\ x^2 & x
\end{array}\right)$
& $ 1.41977\times 10^{-10}$ \\
\hline
$[0.04,0.05]$ &
$\left(\begin{array}{cc}
0 & 1 + 1.00001\,x + 0.499518\,x^2 + 0.174063\,x^3 \\ x^2 & x
\end{array}\right)$
& $1.44084\times 10^{-10}$ \\
\hline
$[0.05,0.06]$ &
$\left(\begin{array}{cc}
0 & 1 + 1.00003\,x + 0.499173\,x^2 + 0.176362\,x^3 \\ x^2 & x
\end{array}\right)$
& $1.44084\times 10^{-10}$ \\
\hline
$[0.06,0.07]$ &
$\left(\begin{array}{cc}
0 & 0.999999 + 1.00004\,x + 0.498952\,x^2 + 0.177587\,x^3 \\ x^2 & x
\end{array}\right)$
& $1.46223\times 10^{-10}$ \\
\hline
$[0.07,0.08]$ &
$\left(\begin{array}{cc}
0 & 0.999998 + 1.00008\,x + 0.498463\,x^2 + 0.179918\,x^3 \\ x^2 & x
\end{array}\right)$
& $1.46223\times 10^{-10}$ \\
\hline
$[0.08,0.09]$ &
$\left(\begin{array}{cc}
 0 & 0.999998 + 1.0001\,x + 0.49816\,x^2 + 0.181181\,x^3 \\ x^2 & x
\end{array}\right)$
& $1.48391\times 10^{-10}$ \\
\hline
$[0.09,0.1]$ &
$\left(\begin{array}{cc}
 0 & 0.999996 + 1.00016\,x + 0.497521\,x^2 + 0.183546\,x^3 \\
   x^2 & x
\end{array}\right)$
& $1.48391\times 10^{-10}$ \\
\hline
\end{tabular}
} \caption{Approximation for Riccati matrix differential
equation~(\ref{paramricc}) in the interval $[0,0.1]$ with step size
$h=0.01$.}
\end{table}
\end{center}

\begin{center}
\begin{figure}
\includegraphics[scale=1.3]{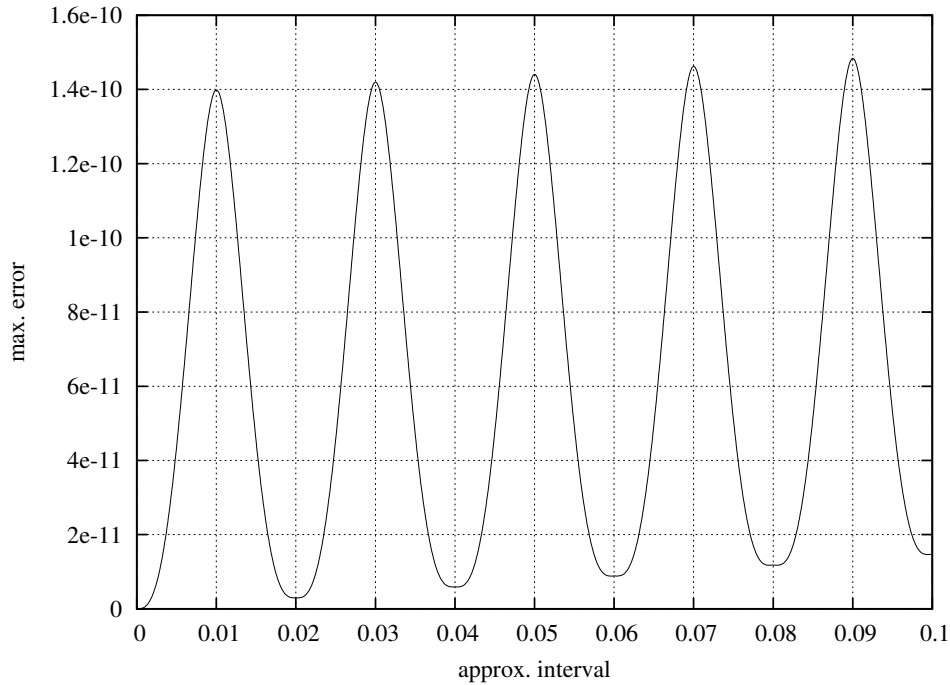}
\caption{Representing the absolute error margins for the Riccati
matrix differential equation~(\ref{paramricc}) in the interval $[0,0.1]$
with step size $h=0.01$.}
\end{figure}
\end{center}

\section{Conclusions}
This article develops a new method for the numerical integration of first-order
matrix differential equations of the non-linear type $Y'(x)=f(x,Y(x)), x\in [a,b]$ using
matrix-cubic splines, and thereby generalizing the approach for the linear case in previous work~\cite{mm}. An important advantage of the proposed method is that the
approximated solution is continuous in the interval under consideration, is
easy to evaluate, and has an error of the order $O(h^4)$. \\ Our method is
well-suited for implementation on numerical and/or symbolical computer systems
({\it Mathematica\/}, {\it Matlab\/}, etc.) as we have shown in Section~3
giving the explicit algorithm.  For a full demonstration of our approach and its
advantages, we conclude with two numerical examples for the Sylvester and
Riccati matrix differential equations.

\end{document}